\documentclass[12pt]{article}

\usepackage{natbib}

\newcommand{\tr}{^{\prime}}

\def\b#1{\mbox{\boldmath $#1$}}    
\def\bl#1{\mbox{\scriptsize \boldmath {$#1$}}} 
\def\cg#1{\mbox{${\cal #1}$}}

\renewcommand{\th}{\theta}
\newcommand{\al}{\alpha}
\newcommand{\be}{\beta}
\newcommand{\de}{\delta}
\newcommand{\la}{\lambda}

\newcommand{\ga}{\gamma}

\newcommand{\diag}{{\rm diag}}    
\newcommand{\pa}{\partial}         
\newcommand{\bu}{{\bar{u}}}         

\def\tb#1{\tilde{\mbox{\boldmath $#1$}}}    

\begin{document}
\title{A note on the application of the Oakes' identity
to obtain the observed information matrix of hidden Markov models}
\author{Francesco Bartolucci, Alessio Farcomeni and Fulvia Pennoni}
\maketitle

\begin{abstract}
We derive the observed information matrix of hidden Markov
models by the application of the Oakes (1999)'s identity. 
The method only requires the first derivative
of the forward-backward recursions of Baum and Welch (1970), instead of
the second derivative of the forward recursion, which is required within
the approach of Lystig and Hughes (2002). The method is illustrated by an example
based on the analysis of a longitudinal dataset which is well known
in sociology.\vspace*{5mm}   

\end{abstract}

\section{Introduction}

Hidden Markov (HM) models have been developed early in the
literature on stochastic processes as
extensions for measurement errors of the standard Markov chain model;
for one of the oldest relevant contributions about these models, see \cite{baum:petr:66}. 
HM models have received much attention in the time-series analysis literature,
due to their wide applicability and easy interpretation
\citep[for an up-to-date review see][]{zucc:macd:09}. These models are also finding an increasing 
popularity for the analysis of longitudinal data \citep[see][]{bart:farc:penn:10}.

The main tool for maximum likelihood (ML) estimation of the parameters
of an HM model is the Expectation-Maximization (EM) algorithm, which is
based on certain forward-backward recursions. This algorithm and these recursions were
developed by Baum and colleagues in a series of papers specifically
for HM models \citep{baum:petr:66,baum:egon:67,baum:et:al:70}. Then, the
EM algorithm was put in a more general context in the widely cited paper of
\cite{demp:lair:rubi:77}.

A drawback of the afore mentioned algorithm is that it does not provide, as a
by-result, the standard errors for the parameter estimates. This is
because it uses neither the observed nor the expected
information matrix, which are suitable transformations of the second
derivative matrix of the model log-likelihood. From the
inverse of these matrices, we obtain standard errors for the
parameter estimates. Then, from the output of the EM algorithm, we
have not an obvious method for assessing the precision of these maximum
likelihood estimates. The information matrix is also important to check
the local identifiability of the model through its rank; see \cite{mchu:56}
and \cite{good:74} among others.

Computing the information matrix (observed or expected) of a 
latent variable model, as an HM model, is considered a difficult task.
Several methods have been proposed to overcome this difficulty;
for a review see \cite{lyst:hugh:02} and \cite{mcla:kris:08}. One of the
more interesting solutions was proposed by \cite{loui:82}.
This solution is based on the missing
information principle as defined by \cite{orch:wood:72}. 
According to this principle, the observed information matrix can be expressed as
the difference between two matrices corresponding to the
complete information, which we would be able to compute if we knew the latent
states, and the missing information due to the unobserved
variables. However, this correction term is difficult in general to compute; see
\cite{oake:99} for further comments and \cite{Turn:Came:Thom:hidd:1998}
for related techniques.

\cite{oake:99} presented an alternative approach, with respect to that
of \cite{loui:82}, to compute the observed information matrix of a latent variable model.
In particular, he derived an explicit formula 
for the second derivative matrix of the model log-likelihood 
which involves the first derivative of the conditional expectation
of the score of the complete data log-likelihood, given observed data.

Specifically for HM models, \cite{lyst:hugh:02} proposed a method
for exactly computing the observed information matrix based on the second
derivative of the forward recursion of \cite{baum:et:al:70} which is used to compute the
model log-likelihood; for a similar method see \cite{bart:06}. 
The method of \cite{lyst:hugh:02} has become rather popular in the HM literature.
Among the methods related to the EM algorithm, we also mention that proposed by
\cite{bart:farc:09} which is very simple to implement and requires a small
extra code over that required for the ML estimation.
However, since it is based on the numerical derivative of the score, the
obtained information matrix may be considered an approximation of the true one.
Also note that, in order to obtain standard errors for the parameter estimates,
we can alternatively use a parametric bootstrap method \citep{efro:tibs:93},
as described in \cite{zucc:macd:09}. Even if the standard errors obtained in this
way may be more reliable with respect to those based on the information matrix,
the method may be computationally costly and, in any case, does not allow us
to check for local identifiability in an obvious way.

In this paper, we show how to apply the \cite{oake:99} identity to
obtain the observed information matrix of an HM model. As we will show, 
the proposed method only requires the first derivative of the
forward-backward recursions of \cite{baum:et:al:70}, whereas the method
of \cite{lyst:hugh:02} requires the second derivative of the forward recursion.
On the other hand, the proposed method is superior to that of \cite{bart:farc:09}
since it allows us to exactly compute the observed information matrix. To the
best of our knowledge, an implementation of the \cite{oake:99}'s identity for
HM models, as the one we propose here, is not available in the literature. 

The proposed approach is illustrated through an application based on a well-known longitudinal
dataset. For the specific HM model used in this application, we make available
some {\tt R} functions\footnote{through a website to be indicated later}
to compute the information matrix and then obtaining the standard errors for the parameter
estimates.

In the following, we first briefly review the EM algorithm and the \cite{oake:99}'s identity
in their general versions. In Section \ref{sec:obs} we propose an implementation
of this identity for HM models on the basis of a suitable reparametrization.
Then, in Section \ref{sec:ex} we describe the application
of the proposed method in connection with the analysis of the dataset mentioned above.

\section{Preliminaries}
We give in this section the necessary background about the EM
algorithm and the \cite{oake:99}'s identity in general; then we
recall some important features about HM models.
\subsection{EM algorithm and observed information matrix}
The EM algorithm \citep{demp:lair:rubi:77} is an iterative algorithm
for finding the ML estimator of models with missing
variables and has a special role in the literature on latent variable models.

With reference to an observed sample, let $\ell(\b\th)$ denote the 
log-likelihood of the latent variable model of interest,
where $\b\th$ is the vector of parameters. As it is well known,
the EM algorithm is based on
the {\em complete data log-likelihood}, denoted as $\ell^*(\b\th)$, 
which is the log-likelihood that we could compute if we knew the 
value of the latent variables for each every sample unit. In particular, 
to maximize $\ell(\b\th)$, the algorithm alternates the following steps 
until convergence:
\begin{itemize}
\item {\bf E-step}: compute the conditional expected value of the complete data
log-likelihood given the current estimate of $\b\th$, denoted by
$\bar{\b\th}$, and the observed data. This expected value is denoted by
$Q(\b\th|\bar{\b\th})$;
\item {\bf M-step}: maximize $Q(\b\th|\bar{\b\th})$ with respect to $\b\th$.
\end{itemize}

We now consider the score and the observed information
matrix corresponding to the model log-likelihood $\ell(\b\th)$.
These are defined, respectively, as
\[
\b s(\b\th) = \frac{\pa\ell(\b\th)}{\pa\b\th}\quad\mbox{and}\quad
\b J(\b\th) = -\frac{\pa^2\ell(\b\th)}{\pa\b\th\pa\b\th\tr}.
\]

It may be simply proved that
\[
\b s(\b\th) = \left.\frac{\pa Q(\b\th|\bar{\b\th})}{\pa\b\th}\right|_{\bar{\bl\th}=\bl\th}.
\]
Consequently, the \cite{oake:99}'s identity states that:
\begin{equation}
\b J(\b\th) = -\left\{
\left.\frac{\pa^2Q(\b\th|\bar{\b\th})}{\pa\b\th\pa\b\th\tr}\right|_{\bar{\bl\th}=\bl\th}+
\left.\frac{\pa^2Q(\b\th|\bar{\b\th})}{\pa\bar{\b\th}\pa\b\th\tr}\right|_{\bar{\bl\th}=\bl\th}
\right\}.\label{eq:oakes}
\end{equation}
This identity then involves two components. The first component is 
the second derivative of the conditional expected value of the complete-data
log-likelihood given the observed data. This component is simple to obtain
from the EM algorithm. The second component involved in (\ref{eq:oakes})
is the first derivative of the score for the same expected log-likelihood
with respect to the current value of the parameters.
\subsection{Hidden Markov models}
Consider now a sequence of $T$ response variables $Y^{(1)},\ldots,Y^{(T)}$,
which are collected in the vector $\b Y$. These response variables may be continuous
or categorical and we may even observe a vector of multivariate outcomes at each $t$.
In the following, we briefly review the assumptions of an HM model for these data
and then how to apply the EM algorithm for ML estimation of the resulting model.
\subsubsection{Assumptions}\label{sec:assumptions}
An HM model relies on the following basic assumptions:
\begin{itemize}
\item the response variables $Y^{(1)},\ldots,Y^{(T)}$ are conditionally
independent given a sequence of unobserved variables $U^{(1)},\ldots,U^{(T)}$ giving rise to
a latent process vector denoted by $\b U$;
\item every response variable $Y^{(t)}$, $t=1,\ldots,T$, depends on
the latent process $\b U$ only through $U^{(t)}$;
\item the latent process $\b U$ follows a Markov chain
with $k$ states labelled from $1$ to $k$.
\end{itemize}
We consider in particular HM models in which:
\begin{itemize}
\item the conditional distribution of $Y^{(t)}$ given $U^{(t)}$ is time-homogenous;
\item the latent Markov chain is of first-order and time-homogeneous.
\end{itemize}

Parameters of the model are then
the initial probabilities of the latent process,
denoted $\la_u=f_{U^{(1)}}(u)$ with $u=1,\ldots,k$ and the transition probabilities
$\pi_{u|\bu}=f_{U^{(t)}|U^{(t-1)}}(u|\bu)$ with $t=2,\ldots,T$ and
$\bu,u=1,\ldots,k$. The initial probabilities are collected in the
$k$-dimensional column vector $\b\la$
and the transition probabilities are collected in the $k\times k$ matrix $\b\Pi$,
with each row denoted $\b\pi_\bu\tr$, where $\b\pi_\bu=
(\pi_{1|\bu},\ldots,\pi_{k|\bu})\tr$.

In the above
expressions, $f_{U^{(t)}}(u)$ denotes the probability mass function of the distribution
of $U^{(t)}$, whereas $f_{U^{(t)}|U^{(t-1)}}(u|\bu)$ denotes the probability mass function
of $U^{(t)}$ given $U^{(t-1)}$. A similar convention will be used
to denote density functions.

Furthermore, when the response variables are categorical
with a reduced number of categories
(labelled from 1 to $c$), we introduce the additional notation
$\phi_{y|u}=f_{Y^{(t)}|U^{(t)}}(y|u)$ with $u=1,\ldots,k$ and $y=0,\ldots,c-1$.
The probabilities are collected in the $c\times k$ matrix $\b\Phi$ which is made
of the column vectors $\b\phi_u$, with $\b\phi=(\phi_{1|u},\ldots,\phi_{c|u})\tr$.
\subsection{Application of the EM algorithm}
It is well known that the above model may be estimated by an EM algorithm 
formulated as in \cite{baum:et:al:70};
see also \cite{bart:farc:penn:10} and \cite{zucc:macd:09}.

Suppose that we observe $n \geq 1$ independent
realizations of $\b Y$, denoted by $\b y_1,\ldots\b y_n$,
with every $\b y_i$
having elements $y_i^{(t)}$, $t=1,\ldots,T$. Note that, in the case
of time-series data, we can only observe a single realization of $\b Y$ and then $n=1$;
in this case, $T$ is typically large.
On the other hand, in the case of longitudinal data,
$n$ is often large as compared to $T$. Our results apply invariably and the model
log-likelihood may be expressed as
\[
\ell(\b\eta) = \sum_i \log f_{\bl Y}(\b y_i)=\sum_{\bl y}n_{\bl y}\log f_{\bl Y}(\b y),
\]
where $\b\eta$ is a vector containing all the parameters in $\b\Phi$, $\b\pi$, and $\b\Pi$,
$f_{\bl Y}(\b y)$ is the probability mass function of $\b Y$ seen
as a function of $\b\eta$. This function can be computed by a forward recursion which
is described in Appendix 1. Moreover, $n_{\bl y}$ is frequency of the response configuration
$\b y=(y^{(1)},\ldots,y^{(T)})\tr$
and the sum $\sum_{\bl y}$ is extended to all response configurations observed
at least once.

We now specialize the EM algorithm for the case of categorical outcomes
mentioned at the end of the previous section.
Let $a^{(t)}_{uy}$, with $t=1,\ldots,T$, $u=1,\ldots,k$, $y=0,\ldots,c-1$,
denote the frequency of $U^{(t)}=u$ and $Y^{(t)}=y$, let $b^{(t)}_u$, with
$t=1,\ldots,T$, $u=1,\ldots,k$, denote the frequency of $U^{(t)}=u$,
and let $c^{(t)}_{\bu u}$, with $t=2,\ldots,T$, $\bu,u=1,\ldots,k$,
denote the joint frequency of the latent states $U^{(t-1)}=\bu$ and $U^{(t)}=u$.
Every E-step of the EM algorithm consists of computing the conditional expected value of
these frequencies given the observed data and the current value of the parameter vector
denoted by $\bar{\b\eta}$, that is
\begin{eqnarray}
\hat{a}^{(t)}_{uy}&=&\sum_i f_{U^{(t)}|\bl Y}(u|\b y_i)I(y_i^{(t)}=y)=
\sum_{\bl y}n_{\bl y}f_{U^{(t)}|\bl Y}(u|\b y)I(y^{(t)}=y),\label{eq:ppost1}\\
\hat{b}^{(t)}_{u}&=&\sum_i f_{U^{(t)}|\bl Y}(u|\b y_i)=
\sum_{\bl y}n_{\bl y}f_{U^{(t)}|\bl Y}(u|\b y),\label{eq:ppost2}\\
\hat{c}^{(t)}_{\bu u}&=&\sum_i f_{U^{(t-1)},U^{(t)}|\bl Y}(y|\b y_i)=
\sum_{\bl y}n_{\bl y}f_{U^{(t-1)},U^{(t)}|\bl Y}(u|\b y),\label{eq:ppost3}
\end{eqnarray}
where $1(\cdot)$ is the indicator function equal to 1 if its argument is true.
These expected values involve posterior probabilities that may be computed
by recursions illustrated in Appendix 1.

Then, the M-step consists of maximizing the conditional expected value,
given the observed data and $\bar{\b\eta}$, of the complete
data log-likelihood, which may be decomposed as
\[
Q(\b\eta|\bar{\b\eta}) = Q_1(\b\Phi|\bar{\b\eta})+Q_2(\b\pi|\bar{\b\eta})+Q_3(\b\Pi|\bar{\b\eta}),
\]
with
\begin{eqnarray*}
Q_1(\b\Phi|\bar{\b\eta})&=&\sum_y \sum_t\sum_u \hat{a}_{uy}^{(t)}\log\phi_{y|u},\\
Q_2(\b\pi|\bar{\b\eta})&=&\sum_u \hat{b}_u^{(1)}\log \la_u,\\
Q_3(\b\Pi|\bar{\b\eta})&=&\sum_{t>1}\sum_{\bu}\sum_u \hat{c}^{(t)}_{\bu u}\log \pi_{u|\bu}.
\end{eqnarray*}
Explicit expressions are available to maximize separately each of
these expressions.
In fact, at every $M$-step $\phi_{y|u}$ is set proportional to $\sum_t
\hat{a}_{uy}^{(t)}$, $\la_u$ to $\hat{b}_u^{(1)}$, and $\pi_{u|\bu}$ to
$\sum_{t>1} \hat{c}^{(t)}_{\bu u}$; see \cite{bart:farc:penn:10}
for a more detailed description.
\section{Observed information matrix for HM models}
\label{sec:obs}
First of all, we consider a reparametrization of the model such that the new
parameter vector, denoted by $\b\th$, is variation independent and is contained in $\cg R^s$
for a suitable $s$. Then, we show how to implement the \cite{oake:99}'s identity
by exploiting this reparametrization.
\subsection{Reparametrization of the model}
The conditional response probabilities are reparametrized through $c-1$ logits
referred to the first category, that is
\[
\al_{y|u} = \log\frac{\phi_{y+1|u}}{\phi_{1|u}},\quad u=1,\ldots,k,\: y=1,\ldots,c-1,
\]
which are included in the $(c-1)$-dimensional column vectors $\b\al_u$; moreover, by
$\b\al$ we denote the vector made of the subvectors $\b\al_1,\ldots,\b\al_k$.
It is worth noting that the choice of
the baseline category is irrelevant for the inference and shall be guided
only by interpretability reasons. The initial probabilities are transformed similarly by
the logits
\[
\be_u = \log\frac{\la_{u+1}}{\la_1},\quad u=1,\ldots,k-1,
\]
which are collected in the $(k-1)$-dimensional column vector $\b\be$. Finally, the
transition probabilities are parametrized through logits referred to the diagonal element,
that is
\[
\ga_{\bu u}=\log\frac{\pi_{u|\bu}}{\pi_{\bu|\bu}},\quad \bu,u=1,\ldots,k,\:u\neq\bu,
\]
which are collected in the $(k-1)$-dimensional 
vectors $\b\ga_\bu$ for $u=1,\ldots,k$; we also denote by
$\b\ga$ the overall vectors made of the subvectors $\b\ga_1,\ldots,\b\ga_k$.

It is convenient to express the above vectors of logits in matrix notation. In particular,
we can easily show that $\b\al$ may be obtained by stacking the vectors
\[
\b\al_u=\b A\log\b\phi_u,\quad u=1,\ldots,k,
\]
where $\b A=\pmatrix{-\b 1_{c-1} & \b I_{c-1}}$, with $\b 1_h$ denoting a column vector of
$h$ ones and $\b I_h$ an identity matrix of the same dimension. The inverse transformation
is
\begin{equation}
\b\phi_u=[\b 1_k\tr\exp(\tb A\b\al_u)]^{-1}\exp(\tb A\b\al_u),\quad
\tb A = \pmatrix{\b 0_{c-1}\tr \cr \b I_{c-1}},\label{eq:inverse_phi}
\end{equation}
where $\b 0_h$ is column vector of $h$ zeros. Similarly, we have that
\[
\b\be=\b B\log\b\la,
\]
with $\b B=\pmatrix{-\b 1_{k-1} & \b I_{k-1}}$; the inverse transformation of
the last expression is
defined as in (\ref{eq:inverse_phi}) on the basis of the matrix $\tilde{\b B}$ defined
in a similar way as $\tilde{\b A}$. Finally, the vector $\b\de$ is made of the subvectors
$\b\de_\bu$, $\bu=1,\ldots,k$, with
\[
\b\de_\bu=\b C_\bu\log\b\pi_\bu,
\]
where
\[
\b C_\bu = \pmatrix{\b I_{\bu-1} & -\b 1_{\bu-1} & \b O_{\bu-1,k-\bu}\cr
                    \b O_{k-\bu,\bu-1} & -\b 1_{k-\bu} & \b I_{k-\bu}},
\]
with $\b O_{hj}$ denoting an $h\times j$ matrix of zeros. The inverse
transformation, to obtain $\b\pi_\bu$ from $\b\de_\bu$, is as in (\ref{eq:inverse_phi}),
with $\tilde{\b A}$ substituted by
\[
\tb C_\bu = \pmatrix{\b I_{\bu-1} & \b O_{\bu-1,k-\bu}\cr
                      \b 0_{\bu-1}\tr & \b 0_{k-\bu}\tr\cr
                      \b O_{k-\bu,\bu-1} & \b I_{k-\bu}}.
\]

The new vector of parameters $\b\th$ is obtained by stacking the single parameters
vectors, that is $\b\th=(\b\al\tr,\b\be\tr,\b\ga\tr)\tr$. Obviously, provided that
all probabilities $\pi_{y|u}$, $\la_u$, and $\pi_{u|\bu}$ are strictly positive,
$\b\th\in\cg R^s$, with $s=(c-1)k + k-1 + k(k-1)$,
and is a one-to-one transformation of the original parameter vector $\b\eta$, which
instead belongs to a more complex space.
\subsection{Computing the observed information matrix}
First of all, adopting the above reparametrization, the expected
value of the complete data log-likelihood may be expressed as
\[
Q(\b\th|\bar{\b\th}) =
Q_1(\b\al|\bar{\b\th})+Q_2(\b\be|\bar{\b\th})+Q_3(\b\ga|\bar{\b\th}),
\]
where, using the matrix notation, we have
\begin{eqnarray*}
Q_1(\b\al|\bar{\b\th})&=&\sum_u\hat{\b a}_u\tr\log\b\phi_u,\\
Q_2(\b\be|\bar{\b\th})&=&(\hat{\b b}^{(1)})\tr\log\b\la,\\
Q_3(\b\ga|\bar{\b\th})&=&\sum_\bu\hat{\b c}_\bu\tr\log\b\pi_\bu,
\end{eqnarray*}
with $\hat{\b a}_u$ denoting a column vector with elements $\sum_t\hat{a}_{uy}^{(t)}$,
$y=0,\ldots,c-1$, $\hat{\b b}^{(1)}$ denoting a column vector with elements
$\hat{b}^{(1)}_u$, $u=1,\ldots,k$, and $\hat{\b c}_\bu$ denoting a vector with elements
$\sum_{t>1}c_{\bu u}^{(t)}$, $u=1,\ldots,k$. Consequently, by applying standard rules
about log-linear models, we have the following score vectors for the complete-data
log-likelihood:
\begin{eqnarray}
\frac{\pa Q_1(\b\al|\bar{\b\th})}{\pa\b\al}&=&
\sum_u\tb A\tr(\hat{\b a}_u-\hat{b}_u^{(t)}\b\phi_u),\label{eq:sc1}\\
\frac{\pa Q_2(\b\be|\bar{\b\th})}{\pa\b\be}&=&
\tb B\tr(\hat{\b b}^{(1)}-n\b\la),\label{eq:sc2}\\
\frac{\pa Q_3(\b\ga|\bar{\b\th})}{\pa\b\ga}&=&
\sum_\bu\tb C_\bu\tr(\hat{\b c}_\bu-\hat{b}_\bu^{(+)}\b\pi_\bu)\label{eq:sc3}.
\end{eqnarray}
where $\b\Omega_{\bl\phi_u} = \diag(\b\phi_u)-\b\phi_u\b\phi_u\tr$,
$\b\Omega_{\bl\la}$ and $\b\Omega_{\bl\pi_\bu}$ are defined in a similar way,
and $\hat{b}_\bu^{(+)}=\sum_{t>1}\hat{b}_\bu^{(t-1)}$.
Similar, we have the second derivative matrices:
\begin{eqnarray*}
\frac{\pa^2Q_1(\b\al|\bar{\b\th})}{\pa\b\al\pa\b\al}&=&
-\sum_u\hat{b}_u^{(t)}\tb A\tr\b\Omega_{\bl\phi_u}\tb A,\\
\frac{\pa^2Q_2(\b\be|\bar{\b\th})}{\pa\b\be\pa\b\be}&=&
-n\tb B\tr\b\Omega_{\bl\la}\tb B,\\
\frac{\pa^2Q_3(\b\ga|\bar{\b\th})}{\pa\b\ga\pa\b\ga}&=&
-\sum_\bu b^{(+)}_\bu\tb C_\bu\tr\b\Omega_{\bl\pi_\bu}\tb C_\bu.
\end{eqnarray*}

It is straightforward to see that the second derivative in (\ref{eq:oakes}) is 
a block-diagonal matrix, with blocks corresponding to above three derivatives,
that is
\[
\frac{\pa^2Q(\b\th|\bar{\b\th})}{\pa\b\th\pa\b\th\tr}=
\diag\left(\frac{\pa^2Q_1(\b\al|\bar{\b\th})}{\pa\b\al\pa\b\al},
\frac{\pa^2Q_2(\b\be|\bar{\b\th})}{\pa\b\be\pa\b\be},
\frac{\pa^2Q_3(\b\ga|\bar{\b\th})}{\pa\b\ga\pa\b\ga}
\right).
\]
Moreover, in order to compute the second component in (\ref{eq:oakes}) we need
the first derivatives of the expected frequencies in (\ref{eq:sc1}), (\ref{eq:sc2}), 
and (\ref{eq:sc3}) with respect to $\bar{\b\th}$.
More precisely, we have
\[
\frac{\pa^2Q(\b\th|\bar{\b\th})}{\pa\bar{\b\th}\pa\b\th\tr}=
\left(
\frac{\pa^2 Q_1(\b\al|\bar{\b\th})}{\pa\bar{\b\th}\pa\b\al\tr},
\frac{\pa^2 Q_2(\b\be|\bar{\b\th})}{\pa\bar{\b\th}\pa\b\be\tr},
\frac{\pa^2 Q_3(\b\ga|\bar{\b\th})}{\pa\bar{\b\th}\pa\b\ga\tr}
\right),
\]
where
\begin{eqnarray}
\frac{\pa^2 Q_1(\b\al|\bar{\b\th})}{\pa\bar{\b\th}\pa\b\al\tr}&=&
\sum_u\left(
\frac{\pa\hat{\b a}_u\tr}{\pa\bar{\b\th}}-
\frac{\pa\hat{b}_u^{(t)}}{\pa\bar{\b\th}}\b\phi_u\tr
\right)\tb A,\label{eq:dsc1}\\
\frac{\pa^2 Q_2(\b\be|\bar{\b\th})}{\pa\bar{\b\th}\pa\b\be\tr}&=&
\frac{\pa(\hat{\b b}^{(1)})\tr}{\pa\bar{\b\th}}\tb B,\label{eq:dsc2}\\
\frac{\pa Q_3(\b\ga|\bar{\b\th})}{\pa\b\ga}&=&
\sum_\bu\left(\frac{\pa\hat{\b c}_\bu\tr}{\pa\bar{\b\th}}-
\frac{\pa\hat{b}_\bu^{(+)}}{\pa\bar{\b\th}}\b\pi_\bu\tr\right)
\tb C_\bu\label{eq:dsc3}.
\end{eqnarray}
How to compute the first derivatives of the above expected values with
respect to $\bar{\b\th}$ is shown in Appendix 2.

Once the observed information at the ML estimate of
$\b\th$ has been obtained through (\ref{eq:oakes}) exploiting the
above results, on the basis of this matrix we can obtain the standard errors
and check identifiability in the usual way. In particular, the standard
errors are obtained by computing the square root of the elements in the
main diagonal of $\b J(\hat{\b\th})^{-1}$. Then, local identifiability is 
checked through the rank of $\b J(\hat{\b\th})$; nevertheless, that
this matrix is of full rank is required in order to compute its inverse.

Note that the standard errors obtained as above are referred to the ML estimate of
the parameter vector $\b\th$. However, we can simply express the standard
errors for the corresponding estimate of the initial parameter vector $\b\eta$
by the delta method. In particular, we first compute
\[
\left(\left.\frac{\pa\b\th\tr}{\pa\b\eta}
\right|_{\bl\eta=\hat{\bl\eta}}\right)\b J(\hat{\b\th})^{-1}
\left(\left.\frac{\pa\b\th}{\pa\b\eta\tr}
\right|_{\bl\eta=\hat{\bl\eta}}\right)
\]
to estimate the variance-covariance matrix of $\hat{\b\eta}$ and then we
obtain the corresponding standard errors as the square root of the elements
in the main diagonal of this matrix. In particular, the derivative matrix
of $\b\th$ with respect of $\b\eta$ may be simply constructed as a block diagonal
matrix with blocks corresponding to the derivative of $\b\al_u$ with respect to
every $\b\phi_u\tr$, $u=1,\ldots,k$, to the derivative of $\b\be$ with respect to
$\b\pi\tr$, and to the derivative of $\b\pi_\bu$ with respect to $\b\de_\bu$,
$\bu = 1,\ldots,k$. For instance, we have
\[
\frac{\pa\b\al_u\tr}{\pa\b\phi}=\b\Omega_{\bl\phi_u}\tb A
\]
and in similar way we can compute the other other blocks.

Finally, it is important to consider that the method described above may be simply
adapted to more sophisticated HM models in which, for instance, the transition
probabilities are time-heterogeneous, the distribution of the response variables
given the latent state is assumed to belong to a certain parametric family, and/or
covariates are included in the model; see \cite{bart:farc:penn:10}. 
However, we prefer to focus on a specific, but important,
HM model in order to make the description of the proposed methods simpler to understand.
\section{Example} \label{sec:ex}
In order to illustrate the proposed approach, we analyze 
a well-known dataset  based on 5 annual waves of the National Youth Survey
\citep{ell:huiz:mena:89}. The dataset concerns 237 individuals who were aged
13 years in 1976. The {\em use of marijuana} was measured by an
ordinal response variable for each wave, having the following three
categories: ``never in the past year" (coded as 1);``no more than
once in a month in the past year" (coded as 2); ``once a month in
the past year" (coded as 3). Such data have been also used for
empirical illustrations by \cite{lang:macd:smit:99},
\cite{verm:hage:04}, and \cite{bart:06}.

With $k =2$ we obtain the estimates of the
conditional response probabilities displayed in Table \ref{tab:cond2esti}.
The table also reports the standard errors obtained with the
proposed method and those obtained using the parametric bootstrap
with a number of sample repetitions equal to 1000. Moreover, in
Tables \ref{tab:ini2esti} and \ref{tab:trans2esti} we show the estimates of
the initial probabilities and of the transition probabilities
respectively, together with the corresponding standard errors.

\begin{table}[ht]
\begin{center}
\begin{tabular}{rccccccccc}
  \hline
 && \multicolumn2c{estimate} && \multicolumn2c{s.e.} && \multicolumn2c{boot.s.e.}\\
 \cline{3-4}\cline{6-7}\cline{9-10}
$y$ && $u=1$ & $u=2$ && $u=1$ & $u=2$ && $u=1$ & $u=2$ \\
  \hline
  1 && 0.9552 & 0.0791 && 0.0137 & 0.0338 && 0.0096 & 0.0315 \\
  2 && 0.0437 & 0.4623 && 0.0131 & 0.0339 && 0.0090 & 0.0338 \\
  3 && 0.0011 & 0.4586 && 0.0024 & 0.0398 && 0.0024 & 0.0358 \\
   \hline
\end{tabular}
\end{center}
\caption{\em Estimates of the parameters $\phi_{y|u}$ and corresponding
standard errors obtained by the proposed method (s.e.) and a parametric bootstrap
method based on 1,000 samples (boot.s.e.).} \label{tab:cond2esti}
\end{table}

\begin{table}[ht]
\begin{center}
\begin{tabular}{rcccccc}
  \hline
 $u$ && est. && s.e. && boot.s.e. \\
  \hline
  1 && 0.9466 && 0.0178 && 0.0166 \\
  2 && 0.0534 && 0.0178 && 0.0166 \\
   \hline
\end{tabular}
\end{center}
\caption{\em Estimates of the parameters $\la_u$ and corresponding
standard errors obtained by the proposed (s.e.) method and a parametric bootstrap
method based on 1,000 samples (boot.s.e.).} \label{tab:ini2esti}
\end{table}

\begin{table}[!ht]
\begin{center}
\begin{tabular}{rccccccccc}
  \hline
 && \multicolumn2c{est.} && \multicolumn2c{se.} && \multicolumn2c{se.boot.}\\
 \cline{3-4} \cline{6-7} \cline{9-10}
 $\bu$ && $u=1$ & $u=2$ && $u=1$ & $u=2$ && $u=1$ & $u=2$ \\
  \hline
  1 && 0.8774 & 0.1226 && 0.0157 & 0.0157 && 0.0140 & 0.0140 \\
  2 && 0.0319 & 0.9681 && 0.0316 & 0.0316 && 0.0268 & 0.0268 \\
   \hline
\end{tabular}
\end{center}
\caption{\em Estimates of the parameters $\pi_{u|\bu}$ and corresponding
standard errors obtained by the proposed method (s.e.) and a parametric bootrap
method based on 1,000 samples (boot.s.e.).} \label{tab:trans2esti}
\end{table}\newpage\vspace*{2mm}

For this application, through the proposed recursion we easily obtain the
standard errors for the parameter estimates. Moreover, as shown in
Tables \ref{tab:cond2esti}, \ref{tab:ini2esti}, and \ref{tab:trans2esti},
these standard errors are always very close to the corresponding parametric
bootstrap standard errors. This confirms the validity of the proposed method
to compute the observed information matrix.

We also estimated the HM model $k=3$ classes, however
the information matrix $\b J(\hat{\b\th})$ is singular
because one of the transition probabilities becomes equal to 0,
so that we cannot state that this model is locally identifiable.

\section*{Appendix 1: Efficient implementation of recursions}
\subsection*{Manifest distribution of the response variables}
In order to efficiently  compute the manifest probability $f_{\bl
Y}(\b y)$, let $\b q^{(t)}(\b y)$ denote the column vector with
elements $f_{U^{(t)},Y^{(1)},\ldots,Y^{(t)}}(u,y^{(1)},\ldots,y^{(t)})$,
for $u=1,\ldots,k$. Then, we have
\begin{equation}
\b q^{(t)}(\b y) = \left\{\begin{array}{ll}
\diag(\b m_{y^{(1)}})\b\la, & t=1,\\
\diag(\b m_{y^{(t)}})\b\Pi\tr\b q^{(t-1)}(\b y), & t=2,\ldots,T,\\
\end{array}\right.\label{eq:forward_recursion}
\end{equation}
where $\b m_y$ is a $k$ dimensional column vector containing the
probabilities $\phi_{y|u}$, $u=1,\ldots,k$.
At the end of this recursions we obtain $f_{\bl Y}(\b y)$ as $\b
q^{(T)}(\b y)\tr\b 1$, where $\b 1$ denotes a column vector of ones
of suitable dimension. In implementing this recursion, attention
must be payed to the case of large values of $T$ because, as $t$
increases, the probabilities in $\b q^{(t)}(\b y)$ could become
negligible; see \cite{Scott:02} for remedial measures.

In the multivariate case, the same recursion as in (\ref{eq:forward_recursion})
may be used, with $\b m_y$ substituted by the vector $\b m_{\bl y}$ with elements
corresponding the conditional probability of the response vector $\b y$ given
every possible value of the corresponding latent state. For further details
on this, and the following recursion, see \cite{zucc:macd:09} and \cite{bart:farc:penn:10}.
\subsection*{Posterior distribution of the latent variables}
Let $\bar{\b q}^{(t)}(\b y)$ be the column vector with elements
$f_{Y^{(t+1)},\ldots,Y^{(T)}|U^{(t)}}(\bar{u},y^{(t+1)},\ldots,y^{(t)})$,
$\bar{u}=1,\ldots,k$. This vector may computed by the backward recursion
\[
\hspace*{-0.5cm} \bar{\b q}^{(t)}(\b y) = \left\{\begin{array}{ll}
\b 1, & t=T,\\
\b\Pi\diag(\b m_{y^{(t+1)}})\bar{\b q}^{(t+1)}(\b y), &
t=T-1,\ldots,1.
\end{array}\right.
\]
Then, the $k$-dimensional column vector $\b f^{(t)}(\b y)$ with
elements $f_{U^{(t)}|\bl Y}(u|\b y)$, $u=1,\ldots,k$, is obtained as
\begin{equation}
\b f^{(t)}(\b y)=\frac{1}{f_{\bl Y}(\b y)}\diag[\b q^{(t)}(\b
y)]\bar{\b q}^{(t)}(\b y),\quad t=1,\ldots,T.\label{eq:deff}
\end{equation}
Moreover, the $k\times k$ matrix $\b F^{(t)}(\b y)$, with elements
$f_{U^{(t-1)},U^{(t)}|\bl Y}(\bar{u},u|\b y)$ arranged by letting
$\bar{u}$ run by row and $u$ by column, is obtained as
\begin{equation}
\b F^{(t)}(\b y) = \frac{1}{f_{\bl Y}(\b y)}\diag[\b q^{(t-1)}(\b
y)]\b\Pi\diag[\b m_{y^{(t)}}]\diag[\bar{\b q}^{(t)}(\b
y)],\label{eq:defF}
\end{equation}
for $t=2,\ldots,T$.
\section*{Appendix 2: derivative of the expected frequencies}
The derivatives of the expected frequencies in (\ref{eq:dsc1}), (\ref{eq:dsc2}),
and (\ref{eq:dsc3}) may be obtained by substituting in (\ref{eq:ppost1}),
(\ref{eq:ppost2}), and (\ref{eq:ppost3}) every posterior probability with
the corresponding derivative with respect to the parameters of interest.
For instance, from (\ref{eq:ppost1}) we have that the derivative matrix
\[
\frac{\pa\hat{\b a}_u\tr}{\pa\bar{\b\th}}
\]
has the following elements 
\[
\frac{\pa\hat{a}^{(t)}_{uy}}{\pa\bar{\th}_j}=
\sum_{\bl y}n_{\bl y}
\frac{\pa f_{U^{(t)}|\bl Y}(u|\b y)}{\pa\bar{\th}_j}
I(y^{(t)}=y),
\]
for $y=1,\ldots,c$, where $\bar{\th}_j$ is an arbitrary element of $\bar{\b\th}$.

In order to compute the derivative of $f_{U^{(t)}|\bl Y}(u|\b y)$, 
and also that of $f_{U^{(t-1)},U^{(t)}|\bl Y}$,
with respect to every parameter $\bar{\th}_j$ we can proceed as 
in \cite{lyst:hugh:02} and \cite{bart:06}. In particular, let
\begin{equation}
\b q^{(t,j)}(\b y)=\frac{\pa \b q^{(t)}(\b y)}{\pa\bar{\th}_j}\quad \mbox{and}\quad
\bar{\b q}^{(t,j)}(\b y)=\frac{\pa\bar{\b q}^{(t)}(\b y)}{\pa\bar{\th}_j},\label{eq:derq}
\end{equation}
and let $\b\phi_u^{(j)}$, $\b\la^{(j)}$, and $\b\Pi^{(j)}$ be defined in a similar way as the
derivatives of $\b\phi_u$, $\b\la$, and $\b\Pi$ with respect to $\bar{\th}_j$;
in a similar way also define $\b m_y^{(j)}$.
Finally, the vectors in (\ref{eq:derq}) may be obtained by the following recursions:
\[
\b q^{(t)}(\b y) = \left\{\begin{array}{ll}
\diag(\b m_{y^{(1)}}^{(j)})\b\la+\diag(\b m_{y^{(1)}})\b\la^{(j)}, & t=1,\\
\diag(\b m_{y^{(t)}}^{(j)})\b\Pi\tr\b q^{(t-1)}(\b y)+
\diag(\b m_{y^{(t)}})(\b\Pi^{(j)})\tr\b q^{(t-1)}(\b y)+\\
+\diag(\b m_{y^{(t)}})\b\Pi\tr\b q^{(t-1,j)}(\b y), & t=2,\ldots,T,\\
\end{array}\right.
\]
and
\[
\hspace*{-0.5cm} \bar{\b q}^{(t)}(\b y) = \left\{\begin{array}{ll}
\b 0 & t=T,\\
\b\Pi^{(j)}\diag(\b m_{y^{(t+1)}})\bar{\b q}^{(t+1)}(\b y)+
\b\Pi\diag(\b m_{y^{(t+1)}}^{(j)})\bar{\b q}^{(t+1)}(\b y)+\\
\b\Pi\diag(\b m_{y^{(t+1)}})\bar{\b q}^{(t+1,j)}(\b y), &
t=T-1,\ldots,1.
\end{array}\right.
\]

Finally, the first derivative of $f_{\bl Y}(\b y)$ with respect to $\bar{\th}_j$
is obtained as $f_{\bl Y}^{(j)}(\b y)=(\b q^{(T,j)})\tr\b 1$.
In a similar way, considering (\ref{eq:deff}) and (\ref{eq:defF}), we obtain
the vector $\b f^{(t,j)}(\b y)$ and $\b F^{(t,j)}(\b y)$, having elements
corresponding to the derivatives of $f_{U^{(t)}|\bl Y}(u|\b y)$ 
and $f_{U^{(t-1)},U^{(t)}|\bl Y}$ with respect to every parameter $\bar{\th}_j$

\bibliographystyle{apalike}
\bibliography{reference}
\end{document}